\documentclass[12pt,twoside]{article}
\usepackage{amsmath}
\usepackage{amssymb}

\def\,{\mskip 3mu} \def\>{\mskip 4mu plus 2mu minus 4mu} \def\;{\mskip 5mu plus 5mu} \def\!{\mskip-3mu}
\def\dispmuskip{\thinmuskip= 3mu plus 0mu minus 2mu \medmuskip=  4mu plus 2mu minus 2mu \thickmuskip=5mu plus 5mu minus 2mu}
\def\textmuskip{\thinmuskip= 0mu                    \medmuskip=  1mu plus 1mu minus 1mu \thickmuskip=2mu plus 3mu minus 1mu}
\textmuskip
\def\beq{\dispmuskip\begin{equation}}    \def\eeq{\end{equation}\textmuskip}
\def\beqn{\dispmuskip\begin{displaymath}}\def\eeqn{\end{displaymath}\textmuskip}
\def\bqa{\dispmuskip\begin{eqnarray}}    \def\eqa{\end{eqnarray}\textmuskip}
\def\bqan{\dispmuskip\begin{eqnarray*}}  \def\eqan{\end{eqnarray*}\textmuskip}

\newtheorem{theorem}{Theorem}

\def\subsection#1{\vspace{1ex}\noindent{\bf{#1.}}}
\def\paragraph#1{\vspace{1ex}\noindent{\bf{#1.}}}

\def\nq{\hspace{-1em}}

\def\qed{\hspace*{\fill}$\Box\quad$}
\def\odt{{\textstyle{1\over 2}}}
\def\eps{\varepsilon}
\def\p{{\scriptscriptstyle+}}
\def\pp{{\scriptscriptstyle++}}
\def\n{{n}}

\def\u{u} %u or \theta or \vartheta

\def\v{\boldsymbol} %\boldsymbol\frak\Bbb\pmb\text
\def\vt{{\v t}}\def\vu{{\v u}}\def\vpi{{\v\pi}}
\def\pin{{\scriptstyle\Pi}}
\def\Var{{\mbox{Var}}}

\def\qmbox#1{{\quad\mbox{#1}\quad}}
\def\propersubset{\subsetneq}
\def\low{\underline}
\def\up{\overline}
\def\DeltaOX{\Delta\hspace{-8pt}{\scriptscriptstyle^{_\otimes}}}
\def\deltaOX{\Delta\hspace{-6.5pt}{\scriptscriptstyle^{_\otimes}}}
\def\Deltapi{\Delta}
\def\Deltal{\Delta_e}
\def\Deltapl{\Delta'_e}
\def\ots{{\scriptscriptstyle\!\otimes}}
\def\Vol{\mbox{Vol}}
\def\argmax{\mathop{\rm arg\,max}}          % maxarg
\def\argmin{\mathop{\rm arg\,min}}          % minarg
\def\leqsq{\sqsubseteq}
\def\geqsq{\sqsupseteq}
\def\Set#1{{\if#1Q{I\!\!\!#1}\else\if#1Z{Z\!\!\!Z}\else{I\!\!#1}\fi\fi}}
\def\maxo{\max}                        %supremum over open Delta
\def\mino{\min}                        %infimum over open Delta

\begin{document}
%%%%%%%%%%%%%%%%%%%%%%%%%%%%%%%%%%%%%%%%%%%%%%%%%%%%%%%%%%%%%%%%%
%                      T i t l e - P a g e                      %
%%%%%%%%%%%%%%%%%%%%%%%%%%%%%%%%%%%%%%%%%%%%%%%%%%%%%%%%%%%%%%%%%

\title{\vskip -5mm\normalsize\sc Technical Report \hfill IDSIA-03-03
\vskip 2mm\bf\LARGE\hrule height5pt \vskip 3mm
\sc Robust Estimators under the \\ Imprecise Dirichlet Model%\thanks{%
\vskip 6mm \hrule height2pt \vskip 5mm}
\author{{\bf Marcus Hutter}\\[3mm]
\normalsize IDSIA, Galleria 2, CH-6928\ Manno-Lugano, Switzerland\\
\normalsize marcus@idsia.ch \hspace{8.5ex} http://www.idsia.ch/$^{_{_\sim}}\!$marcus}
\date{May 5, 2003}
\maketitle

\begin{abstract}
Walley's Imprecise Dirichlet Model (IDM) for categorical data
overcomes several fundamental problems which other approaches to
uncertainty suffer from. Yet, to be useful in practice, one needs
efficient ways for computing the imprecise=robust sets or
intervals. The main objective of this work is to derive exact,
conservative, and approximate, robust and credible interval
estimates under the IDM for a large class of statistical
estimators, including the entropy and mutual information.
\end{abstract}
\vspace{5cm}

\pagebreak

%%%%%%%%%%%%%%%%%%%%%%%%%%%%%%%%%%%%%%%%%%%%%%%%%%%%%%%%%%%%%%%
\section{Introduction}\label{secInt}
%%%%%%%%%%%%%%%%%%%%%%%%%%%%%%%%%%%%%%%%%%%%%%%%%%%%%%%%%%%%%%%
This work derives interval estimates under the Imprecise Dirichlet
Model (IDM) \cite{Walley:96} for a large class of statistical
estimators. In the IDM one considers an i.i.d.\ process with
unknown chances\footnote{Also called {\em objective} or {\em
aleatory} probabilities.} $\pi_i$ for outcome $i$. The prior
uncertainty about $\vpi$ \footnote{We denote vectors by $\v
x:=(x_1,...,x_d)$ for $\v x\in\{\v n,\vt,\vu,\vpi,...\}$.}is
modeled by a set of Dirichlet priors\footnote{Also
called {\em second order} or {\em subjective} or {\em belief} or
{\em epistemic} probabilities.}
$\{p(\vpi)\propto\prod_i\pi_i^{st_i-1}\,:\,\vt\in\Delta\}$, where%
\footnote{Strictly speaking, $\Delta$ should be the open simplex
\cite{Walley:96}, since $p(\vpi)$ is improper for $\vt$ on the
boundary of $\Delta$. For simplicity we assume that, if necessary,
considered functions of $\vt$ can and are continuously extended
to the boundary of $\Delta$, so that, for instance, minima and
maxima exist. All considerations can straightforwardly, but
cumbersomely, be rewritten in terms of an open simplex. Note that
open/closed $\Delta$ result in open/closed robust intervals, the
difference being numerically/practically irrelevant.
} $\Delta:=\{\vt\,:\,t_i\geq 0,\, \sum_i t_i=1\}$, and $s$ is a
hyper-parameter, typically chosen between 1 and 2. Sets of
probability distributions are often called Imprecise
probabilities, hence the name IDM for this model. We avoid the
term {\em imprecise} and use {\em robust} instead, or capitalize
{\em Imprecise}. IDM overcomes several fundamental problems which
other approaches to uncertainty suffer from \cite{Walley:96}. For
instance, IDM satisfies the representation invariance principle
and the symmetry principle, which are mutually exclusive in a pure
Bayesian treatment with proper prior \cite{Walley:96}. The counts
$n_i$ for $i$ form a minimal sufficient statistic of the data of
size $n=\sum_i n_i$. Statistical estimators $F(\v n)$ usually also
depend on the chosen prior: so a set of priors leads to a set of
estimators $\{F_\vt(\v n)\,:\;\vt\in\Delta\}$. For instance,
the expected chances $E_\vt[\pi_i]={n_i+st_i\over
n+s}=:\u_i(\vt)$ lead to a robust interval estimate $[{n_i\over
n+s},{n_i+s\over n+s}]\ni E_\vt[\pi_i]$. Robust intervals for
the variance $\Var[\pi_i]$ \cite{Walley:96} and for the mean and
variance of linear-combinations $\sum_i\alpha_i\pi_i$ have also
been derived \cite{Bernard:01}. Bayesian estimators (like
expectations) depend on $\vt$ and $\v n$ only through $\vu$
(and $n+s$ which we suppress), i.e.\ $F_\vt(\v n)=F(\vu)$.
The main objective of this work is to derive approximate,
conservative, and exact intervals $[\mino_{\vt\in\Delta}F(\vu),\maxo_{\vt\in\Delta}F(\vu)]$ for general
$F(\vu)$, and for the expected (also called predictive) entropy
and the expected mutual information in particular. These results
are key building blocks for applying IDM. Walley suggests, for
instance, to use $\mino_\vt P_\vt[{\cal F}\geq c]\geq\alpha$
for inference problems and $\mino_\vt E_\vt[{\cal F}]\geq c$
for decision problems \cite{Walley:96}, where $\cal F$ is some
function of $\vpi$. One application is the inference of robust
tree-dependency structures \cite{Zaffalon:01tree,Hutter:03tree},
in which edges are partially ordered based
on Imprecise mutual information.

Section \ref{secIDM} gives a brief introduction to
IDM and describes our problem setup.
In Section \ref{secEEI} we derive exact robust intervals for
concave functions $F$, such as the entropy.
Section \ref{secAEI} derives approximate robust intervals for arbitrary
$F$.
In Section \ref{secEP} we show how bounds of elementary functions
can be used to get bounds for composite function, especially for
sums and products of functions. The results are used in
Section \ref{secIEMI} for deriving robust intervals for the
mutual information.
The issue of how to set up IDM models on
product spaces is discussed in Section \ref{secPS}.
Section \ref{secCI} addresses the problem of how to combine
Bayesian credible intervals with the robust intervals of the IDM.
Conclusions are given in Section \ref{secConc}.

%%%%%%%%%%%%%%%%%%%%%%%%%%%%%%%%%%%%%%%%%%%%%%%%%%%%%%%%%%%%%%%
\section{The Imprecise Dirichlet Model}\label{secIDM}
%%%%%%%%%%%%%%%%%%%%%%%%%%%%%%%%%%%%%%%%%%%%%%%%%%%%%%%%%%%%%%%

%-------------------------------%
\subsection{Random i.i.d.\ processes}
%-------------------------------%
We consider discrete random variables $\imath\in\{1,...,d\}$ and
an i.i.d.\ random process with outcome $i\in\{1,...,d\}$ having
probability $\pi_i$. The chances $\vpi$ form a probability
distribution, i.e.\ $\vpi\in\Deltapi:=\{\v x\in\Set R^d\,:\,x_i\geq
0\,\forall i,\; x_\p=1\}$, where we have used the abbreviation $\v
x=(x_1,...,x_d)$ and $x_\p:=\sum_{i=1}^d x_i$. The likelihood of a
specific data set $\v D$ with $n_i$ observations $i$ and total
sample size $n=n_\p=\sum_i n_i$ is $p(\v
D|\vpi)=\prod_i\pi_i^{n_i}$.
The chances $\pi_i$ are usually unknown and have to be estimated
from the sample frequencies $n_i$. The frequency estimate
${n_i\over n}$ for $\pi_i$ is one possible point estimate.

%-------------------------------%
\subsection{Second order p(oste)rior}
%-------------------------------%
In the Bayesian approach one models the initial uncertainty in
$\vpi$ by a (second order) prior ``belief'' distribution
$p(\vpi)$ with domain $\vpi\in\Deltapi$. The Dirichlet priors
$p(\vpi)\propto\prod_i\pi_i^{n'_i-1}$, where $n'_i$ comprises prior
information, represent a large class of priors. $n'_i$ may be
interpreted as (possibly fractional) virtual number of ``observation''.
High prior belief in $i$ can be modeled by large $n'_i$. It is
convenient to write $n'_i=s\cdot t_i$ with $s:=n'_+$, hence $\vt\in\Delta$. Having no initial bias one should choose a prior in
which all $t_i$ are equal, i.e.\ $t_i={1\over d}\,\forall i$.
Examples for $s$ are $0$ for Haldane's prior \cite{Haldane:48},
$1$ for Perks' prior \cite{Perks:47}, ${d\over 2}$ for Jeffreys'
prior \cite{Jeffreys:46}, and $d$ for Bayes-Laplace's uniform
prior \cite{Gelman:95}. From the prior and the data likelihood one
can determine the posterior $p(\vpi|\v D)=p(\vpi|\v
n)\propto\prod_i\pi_i^{n_i+st_i-1}$.

The posterior $p(\vpi|\v D)$ summarizes all statistical
information available in the data. In general, the posterior is a
very complex object, so we are interested in summaries of this
plethora of information. A possible summary is the expected
value or mean
$
  E_\vt[\pi_i]={n_i+st_i\over n+s}
$
which is often used for estimating $\pi_i$. The accuracy may be
obtained from the covariance of $\vpi$.

Usually one is not only interested in an estimation of the whole
vector $\vpi$, but also in an estimation of scalar functions
${\cal F}:\Deltapi\to I\!\!R$ of $\vpi$, such as the entropy
${\cal H}(\vpi)=-\sum_i\pi_i\log\pi_i$, where $\log$ denotes the
natural logarithm. Since $\cal F$ is itself a random variable we
could determine the posterior distribution $p({\cal F}_0|\v
n)=\int_{\Deltapi}\delta({\cal F}(\vpi)-{\cal F}_0)p(\vpi|\v
n)d\vpi$ of ${\cal F}$, which may further be summarized by the
posterior mean $E_\vt[{\cal F}]=\int_{\Deltapi}{\cal
F}(\vpi)p(\vpi|\v n)d\vpi$ and possibly the posterior variance
$\Var[{\cal F}]$.
A simple, but crude approximation for the mean can be obtained by
exchanging $E$ with ${\cal F}$ (exact only for linear functions):
$E_\vt[{\cal F}(\vpi)]\approx {\cal F}(E_\vt[\vpi])$. The
approximation error is typically of the order $1\over n$.

%-------------------------------%
\subsection{The Imprecise Dirichlet Model}
%-------------------------------%
The classical approach, which consists of selecting a single
prior, suffers from a number of problems. Firstly, choosing for
example a uniform prior $t_i={1\over d}$, the prior depends on the
particular choice of the sampling space.
Secondly, it assumes exact prior knowledge of $p(\vpi)$. The
solution to the second problem is to model our ignorance by
considering sets of priors $p(\vpi)$, often called Imprecise
probabilities. The specific {\em Imprecise Dirichlet Model} (IDM)
\cite{Walley:96} considers the set of {\em all} $\vt\in\Delta$,
i.e. $\{p(\vpi|\v n):\vt\in\Delta\}$ which solves also the first
problem. Walley suggests to fix the hyperparameter $s$ somewhere
in the interval $[1,2]$. A set of priors results in a set of
posteriors, set of expected values, etc. For real-valued
quantities like the expected entropy $E_\vt[{\cal H}]$ the sets
are typically intervals, which we call robust
intervals
\beqn
  E_\vt[{\cal F}] \in [\mino_{\vt\in\Delta}E_\vt[{\cal F}] \,,\,
    \maxo_{\vt\in\Delta}E_\vt[{\cal F}]].
\eeqn

%-------------------------------%
\subsection{Problem setup and notation}
%-------------------------------%
Consider any statistical estimator $F$. $F$ is a function of the
data $\v D$ and the hyperparameters $\vt$.
We define the general correspondence
\beq\label{thtrel}
 \u_i^{\cdots}={n_i+st_i^{\cdots}\over n+s},
 \quad\mbox{where $^{\ldots}$ can be various superscripts}.
\eeq
$F$ can, hence, be rewritten as a function of $\vu$ and $\v D$.
Since we regard $\v D$ as fixed, we suppress this dependence and
simply write $F=F(\vu)$. This is further motivated by the fact
that all Bayesian estimators of functions $\cal F$ of $\vpi$ only
depend on $\vu$ and the sample size $n+s$. It is easy to see
that this holds for the mean, i.e.\ $E_\vt[{\cal
F}]=F(\vu\,;\,n+s)$, and similarly for the variance and all higher
(central) moments. The main focus of this work is to derive exact
and approximate expressions for upper and lower $F$ values
\beqn
 \up{F}:=\maxo_{\vt\in\Delta}F(\vu)
 \qmbox{and}
 \low F:=\mino_{\vt\in\Delta}F(\vu),\qquad
 \up{\low F}:=[\low F,\up F]
\eeqn
$\vt\in\Delta$ $\Leftrightarrow$
$\vu\in\Delta'$, where $\Delta':=\{\vu\,:\,\u_i\geq{n_i\over
n+s},\;\u_\p=1\}$. We define $\vu^{\up F}$ as the
$\vu\in\Delta'$ which maximizes $F$, i.e.\ $\up F=F(\vu^{\up
F})$, and similarly $\vt^{\up F}$ through relation
(\ref{thtrel}). If the maximum of $F$ is assumed in a corner of
$\Delta'$ we denote the index of the corner by $i^{\up F}$, i.e.\
$t_i^{\up F}=\delta_{ii^{\up F}}$, where $\delta_{ij}$ is
Kronecker's delta function. Similarly $\vu^{\low F}$, $\vt^{\low F}$, $i^{\low F}$.

%%%%%%%%%%%%%%%%%%%%%%%%%%%%%%%%%%%%%%%%%%%%%%%%%%%%%%%%%%%%%%%
\section{Exact Robust Intervals for Concave Estimators}\label{secEEI}
%%%%%%%%%%%%%%%%%%%%%%%%%%%%%%%%%%%%%%%%%%%%%%%%%%%%%%%%%%%%%%%

In this section we derive exact expressions for $\up{\low F}$ if
$F:\Delta\to I\!\!R$ is of the form
\beq\label{Fconc}
  F(\vu)=\sum_{i=1}^d f(\u_i)
  \qmbox{and concave} f:[0,1]\to I\!\!R.
\eeq
The expected entropy is such an example (discussed later). Convex
$f$ are treated similarly (or simply take $-f$).

%-------------------------------%
\subsection{The nature of the solution}
%-------------------------------%
The approach to a solution of this problem is motivated as
follows: Due to symmetry and concavity of $F$, the global maximum
is attained at the center $\u_i={1\over d}$ of the probability
simplex $\Delta$, i.e.\ the more uniform $\vu$ is, the larger
$F(\vu)$. The nearer $\vu$ is to a vertex of $\Delta$, i.e.\
the more unbalanced $\vu$ is, the smaller is $F(\vu)$. The
constraints $t_i\geq 0$ restrict $\vu$ to the smaller simplex
\beqn
  \Delta'=\{\vu\,:\,\u_i\geq\u_i^0,\;\u_\p=1\}
  \qmbox{with} \u_i^0:={n_i\over n+s},
\eeqn
which prevents setting $\u_i^{\up
F}={1\over d}$ and $\u_i^{\low F}=\delta_{i1}$.
Nevertheless, the basic idea of choosing $\vu$ as
uniform / as unbalanced as possible still works, as we will see.

%-------------------------------%
\subsection{Greedy $F(\vu)$ minimization}
%-------------------------------%
Consider the following procedure for obtaining $\vu^{\low F}$.
We start with $\vt\equiv\v 0$ (outside the usual domain $\Delta$
of $F$, which can be extended to
$[0,1]^d$ via (\ref{Fconc})) and then gradually increase $\vt$ in
an axis-parallel way until $t_\p=1$. With axis-parallel we mean
that only one component of $\vt$ is increased, which one possibly
changes during the process. The total zigzag curve from $\vt^{start}=\v 0$ to $\vt^{end}$ has length $t_\p^{end}=1$. Since
all possible curves have the same (Manhattan) length 1,
$F(\vu^{end})$ is minimized for the curve which has (on average)
smallest $F$-gradient along its path. A greedy strategy is to
follow the direction $i$ of currently smallest $F$-gradient
${\partial F\over\partial t_i}=f'(\u_i){s\over n+s}$. Since $f'$
is monotone decreasing ($f''<0$), ${\partial F\over\partial t_i}$
is smallest for largest $\u_i$. At $\vt^{start}=\v 0$,
$\u_i={n_i\over n+s}$ is largest for $i=i^{min}:=\arg\max_i n_i$.
Once we start in direction $i^{min}$, $\u_{i^{min}}$ increases
even further whereas all other $\u_i$ ($i\neq i^{min}$) remain
constant. So the moving direction is never changed and finally we
reach a local minimum at $t_i^{end}=\delta_{ii^{min}}$.
In \cite{Hutter:03idmx} we show that this is a global minimum,
i.e.\
\beq\label{Fmin}
  t_i^{\low F}=\delta_{ii^{\low F}} \qmbox{with}
  i^{\low F}:=\arg\max_i n_i.
\eeq

%-------------------------------%
\subsection{Greedy $F(\vu)$ maximization}
%-------------------------------%
Similarly we maximize $F(\vu)$. Now we increase $\vt$ in direction
$i=i_1$ of maximal ${\partial F\over\partial t_i}$, which is
the direction of smallest $\u_i\propto n_i+st_i$.
Again, (only) $\u_{i_1}$ increases, but possibly reaches a value
where it is no longer the smallest one. We stop if it becomes
equal to the second smallest $\u_i$, say $i=i_2$. We now have to
increase $\u_{i_1}$ and $\u_{i_2}$ with same speed (or in an
$\eps$-zigzag fashion) until they become equal to $\u_{i_3}$, etc
or until $\u_\p=1=t_\p$ is reached. Assume the process stops with
direction $i_m$ and minimal $\u$ being $\tilde\u$, i.e.\ finally
$\u_{i_k}=\tilde\u$ for $k\leq m$ and $t_{i_k}=0$ for $k>m$.
From the constraint $1=\u_\p=\sum_{k\leq
m}\u_{i_k}+\sum_{k>m}\u_{i_k} = m\tilde\u+\sum_{k>m}{n_{i_k}\over
n+s}$ we obtain $\tilde\u(m)={1\over m}[1-\sum_{k>m}{n_{i_k}\over n+s}]=
[s+\sum_{k\leq m}n_{i_k}]/[m(n+s)]$. One can show that
$\tilde\u(m)$ has one global minimum (no local ones)
and that the final $m$
is the one which minimizes $\tilde\u$, i.e.\
\beq\label{Fmax}
  \tilde\u = \min_{m\in\{1...d\}} {s+\sum_{k\leq m} n_{i_k}\over
  m(n+s)}, \;\;\mbox{where}\;\; n_{i_1}\leq
  n_{i_2}\leq...\leq\n_{i_d},\quad
  \u_i^{\up F}=\max\{\u_i^0,\tilde\u\}.
\eeq
If there is a unique minimal $n_{i_1}$ with gap
$\geq s$ to the 2nd smallest $n_{i_2}$ (which is quite likely for not
too small $n$), then $m=1$ and the maximum is attained at a
corner of $\Delta$ ($\Delta'$).

\begin{theorem}[Exact extrema for concave functions on simplices]\label{thEEI}
Assume $F:\Delta'\to I\!\!R$ is a concave function of the form
$F(\vu)=\sum_{i=1}^df(\u_i)$.
%, where $\Delta'=\{\vu\,:\,\u_i\geq {n_i\over n+s},\;\u_\p=1\}$.
Then $F$ attains the global maximum $\up F$ at $\vu^{\up F}$
defined in (\ref{Fmax}) and the global minimum $\low F$ at
$\vu^{\low F}$ defined in (\ref{Fmin}).
\end{theorem}

\paragraph{Proof} What remains to be shown is that the solutions
obtained in the last paragraphs by greedy
minimization/maximization of $F(\vu)$ are actually global
minima/maxima. For this assume that $\vt$ is a local minimum of
$F(\vu)$. Let $j:=\arg\max_i\u_i$ (ties broken arbitrarily).
Assume that there is a $k\neq j$ with non-zero $t_k$. Define
$\vt'$ as $t'_i=t_i$ for all $i\neq j,k$, and $t'_j=t_j+\eps$,
$t'_k=t_k-\eps$, for some $0<\eps\leq t_k$. From $\u_k\leq \u_j$
and the concavity of $f$ we get\footnote{Slope
${f(\u+\eps)-f(\u)\over\eps}$ is a decreasing function in $\u$
for any $\eps>0$, since $f$ is concave.}
\bqan
  F(\vt')-F(\vt) &=& [f(\u'_j)+f(\u'_k)]-[f(\u_j)+f(\u_k)]
  \\ &=&
  [f(\u_j\!+\!\sigma\eps)-f(\u_j)]-[f(\u_k)-f(\u_k\!-\!\sigma\eps)] \;<\; 0
\eqan
where $\sigma:={s\over n+s}$. This contradicts the minimality
assumption of $\vt$. Hence, $t_i=0$ for all $i$ except one (namely $j$,
where it must be 1). (Local) minima are attained in the vertices
of $\Delta$. Obviously the global minimum is for $t_i^{\low
F}=\delta_{ii^{\low F}}$ with $i^{\low F}:=\arg\max_i n_i$. This
solution coincides with the greedy solution. Note that the global
minimum may not be unique, but since we are only interest in
the value of $F(\vu^{\low F})$ and not its argument this
degeneracy is of no further significance.

Similarly for the maximum, assume that $\vt$ is a (local)
maximum of $F(\vu)$. Let $j:=\arg\min_i\u_i$ (ties broken
arbitrarily). Assume that there is a $k\neq j$ with
non-zero $t_k$ {\em and} $\u_k>\u_j$. Define $\vt'$ as above with
$0<\eps<\min\{t_k\,,\,t_k-t_j\}$. Concavity of $f$ implies
\beqn
  F(\vt')-F(\vt) =
  [f(\u_j\!+\!\sigma\eps)-f(\u_j)]-[f(\u_k)-f(\u_k\!-\!\sigma\eps)] > 0,
\eeqn
which contradicts the maximality assumption of $\vt$. Hence
$t_i=0$ if $\u_i$ is not minimal ($\tilde\u$). The previous
paragraph constructed the unique solution $\vu^{\up F}$ satisfying
this condition. Since this is the only local maximum it must be
the unique global maximum (contrast this to the minimum case).
\qed

\begin{theorem}[Exact extrema of expected entropy]\label{corEEI}
Let ${\cal H}(\vpi)=-\sum_i \pi_i\log\pi_i$ be the entropy of
$\vpi$ and the uncertainty of $\vpi$ be modeled by the Imprecise
Dirichlet Model. The expected entropy $H(\vu):=E_\vt[{\cal
H}]$ for given hyperparameter $\vt$ and sample $\v n$ is given by
\beq\label{hex}
  H(\vu)=\sum_i h(\u_i) \qmbox{with}
  h(\u)=\u\!\cdot\![\psi(n\!+\!s\!+\!1)-\psi((n\!+\!s)\u\!+\!1)] =
  \u\cdot\!\!\!\!\nq\sum_{k={(n+s)}u+1}^{n+s}\nq\!\! k^{-1}
\eeq
where $\psi(x)=d\,\log\Gamma(x)/dx$ is the logarithmic derivative of
the Gamma function and the last expression is valid for integral $s$
and $(n+s)u$. The lower $\low H$ and upper $\up H$ expected
entropies are assumed at $\vu^{\low H}$ and $\vu^{\up H}$ given in
(\ref{Fmin}) and (\ref{Fmax}) (with $F\leadsto H$, see also
(\ref{thtrel})).
\end{theorem}

A derivation of the exact expression (\ref{hex}) for the
expected entropy can be found in
\cite{Wolpert:95,Hutter:01xentropy}. The only thing to be shown is
that $h$ is concave. This may be done by exploiting special
properties of the digamma function $\psi$ (see
\cite{Abramowitz:74}).

There are fast implementations of $\psi$ and its derivatives and
exact expressions for integer and half-integer arguments

\paragraph{Example} For $d=2$, $n_1=3$, $n_2=6$, $s=1$ we have $n=9$,
$\u_1={3+t_1\over 10}$, $\u_2={6+t_2\over 10}$, $\vt^0=0$,
$\vu^0=({.3\atop .6})$, see (\ref{thtrel}).
From (\ref{Fmin}), $i^{\low H}=2$, $\vt^{\low
H}=({0\atop 1})$, $\vu^{\low H}=({.3\atop .7})$.
From (\ref{Fmax}), $i_1=1$, $i_2=2$, $\tilde\u=\min\{{1+3\over
9+1},{1+3+6\over 2\cdot(9+1)}\}={4\over 10}$, $\vu^{\up
H}=\max\{\vu^0,\tilde\u\}=({.4\atop .6})$ $\Rightarrow$
$\vt^{\up H}=({1\atop 0})$ is in corner.
From (\ref{hex}), %
$h({3\over 10})={2761\over 8400}$,
$h({4\over 10})={2131\over 6300}$,
$h({6\over 10})={1207\over 4200}$,
$h({7\over 10})={847\over 3600}$, hence
$\up{\low H}=[H(\vu^{\low H}),H(\vu^{\up H})]=
 [h({3\over 10})+h({7\over 10})\, ,\, h({4\over 10})+h({6\over 10})]
 =[0.5639...,0.6256...]$, so $\up H-\low H=O({1\over 10})$.

%%%%%%%%%%%%%%%%%%%%%%%%%%%%%%%%%%%%%%%%%%%%%%%%%%%%%%%%%%%%%%%
\section{Approximate Robust Intervals}\label{secAEI}
%%%%%%%%%%%%%%%%%%%%%%%%%%%%%%%%%%%%%%%%%%%%%%%%%%%%%%%%%%%%%%%

%-------------------------------%
%\subsection{Taylor expansion of $F(\vu)$}
%-------------------------------%
In this section we derive approximations for $\up{\low F}$
suitable for arbitrary, twice differentiable functions
$F(\vu)$.
The derived approximations for $\up{\low F}$ will be robust in the
sense of covering set $\up{\low F}$ (for any $n$), and the
approximations will be ``good'' if $n$ is not too small.
In the following, we treat $\sigma:={s\over n+s}$ as a (small)
expansion parameter.
For $\vu,\vu^*\in\Delta'$ we have
\beq\label{dtbnd}
  \u_i-\u_i^* \;=\; \sigma\!\cdot\!(t_i-t_i^*) \qmbox{and}
  |\u_i-\u_i^*| \;=\; \sigma|t_i-t_i^*| \;\leq\; \sigma
  \qmbox{with} \sigma:=\textstyle{s\over n+s}.
\eeq
Hence we may Taylor-expand $F(\vu)$ around $\vu^*$, which
leads to a Taylor series in $\sigma$. This shows that $F$ is
approximately linear in $\vu$ and hence in $\vt$. A linear
function on a simplex assumes its extreme values at the vertices
of the simplex. This has already been encountered in Section
\ref{secEEI}. The consideration above is a simple explanation for
this fact. This also shows that the robust interval $\up{\low F}$
is of size $\up F-\low F=O(\sigma)$.\footnote{$f(\v n,\vt,s)=O(\sigma^k)$ $\;:\Leftrightarrow\;$ $\exists c\,\forall\v
n\in I\!\!N_0^d,\,\vt\in\Delta,\,s>0$ : $|f(\v n,\vt,s)|\leq
c\sigma^k$, where $\sigma={s\over n+s}$.} Any approximation to
$\up{\low F}$ should hence be at least $O(\sigma^2)$. The
expansion of $F$ to $O(\sigma)$ is
\beq\label{Fexpand}
  F(\vu) \;=\; \overbrace{F(\vu^*)}^{F_0=O(1)} +
  \overbrace{\sum_i[\partial_i F(\v{\check\u})](\u_i-\u_i^*)}^{F_R=O(\sigma)}
\eeq
where $\partial_i F(\v{\check\u})$ is the partial derivative
$\partial_i F(\v{\check\u})\over\partial\check\u_i$ of
$F(\v{\check\u})$ w.r.t.\ $\check\u_i$. For suitable
$\v{\check\u}=\v{\check\u}(\vu,\vu^*)\in\Delta'$ this expansion is
exact ($F_R$ is the exact remainder).
Natural points for expansion are $t_i^*={1\over d}$ in the center
of $\Delta$, or possibly also $t_i^*={n_i\over n}=\u_i^*$. See
\cite{Hutter:03idmx}
for such a general expansion. Here, we expand around the improper
point $t_i^*:=t_i^0\equiv 0$, which is outside(!) $\Delta$, since
this makes expressions particularly
simple.$\!\!$\footnote{The order of accuracy
$O(\sigma^2)$ we will encounter is for all choices of $\vu^*$
the same. The concrete numerical errors differ of course. The
choice $\vt^*=\v 0$ can lead to $O(d)$ smaller $F_R$ than the
natural center point $\vt^*={1\over d}$, but is more likely a
factor $O(1)$ larger. The exact numerical values depend on the
structure of $F$.} (\ref{dtbnd}) is still valid in this case,
and $F_R$ is exact for some
$\v{\check\u}$ in
\beqn
  \Deltapl:=\{\vu\,:\,\u_i\geq\u_i^0\,\forall i,\;\u_\p\leq 1\},
  \qmbox{where}
  \u_i^0={n_i\over n+s}.
\eeqn
Note that we keep the exact condition
$\vu\in\Delta'$. $F$ is usually already defined on
$\Deltapl$ or extends from $\Delta'$ to $\Deltapl$
without effort in a natural way (analytical continuation). We
introduce the notation
\beq\label{eqleqsq}
  F\leqsq G \quad:\Leftrightarrow\quad
  F\leq G \qmbox{and} F=G+O(\sigma^2)
\eeq
stating that $G$ is a ``good'' upper bound on $F$.
The following bounds hold for arbitrary differentiable functions.
In order for the bounds to be ``good,$\!$'' $F$ has to be Lipschitz
differentiable in the sense that there exists a constant $c$ such that
\beqn
  |\partial_i F(\vu)|\leq c \qmbox{and}
  |\partial_i F(\vu)-\partial_i F(\vu')|
  \leq c|\vu-\vu'|
\eeqn
\beq\label{Lipschitz}
  \forall\,\vu,\vu'\in\Deltapl
  \qmbox{and} \forall\,1\leq i \leq d.
\eeq
\beqn
  \mbox{If $F$ depends also on $\v n$, e.g.\ via $\sigma$ or
  $\vu^0$, then $c$ shall be independent of them.}
\eeqn
The Lipschitz condition is satisfied, for instance, if the
curvature $\partial^2 F$ is uniformly bounded. This is satisfied
for the expected entropy $H$ (see (\ref{hex})), but violated for
the approximation $E_\vt[{\cal H}]\approx{\cal H}(\vu)$ if $n_i=0$ for
some $i$.

\begin{theorem}[Approximate robust intervals]\label{thARI}
Assume $F:\Deltapl\to I\!\!R$ is a Lipschitz differentiable
function (\ref{Lipschitz}).
Let $[\low F,\up F]$ be the global [minimum,maximum] of $F$ restricted
to $\Delta'$. Then
\beqn
  F(\vu^1) \;\leqsq\; \up F \;\leqsq\; F_0+F_R^{ub}
  \;\;\mbox{where}\;\; F_R^{ub}=\max_i F_{iR}^{ub}
  \;\mbox{and}\;\;
  F_{iR}^{ub} \;=\; \sigma\maxo_{\vu\in\Deltapl}
    [\partial_i F(\vu)]
\eeqn\vspace{-1ex}
\beqn
  F_0+F_R^{lb} \;\leqsq\; \low F \;\leqsq\; F(\vu^2)
  \;\;\mbox{where}\;\; F_R^{lb}=\min_i F_{iR}^{lb}
  \;\;\mbox{and}\;\;
  F_{iR}^{lb} \;=\; \sigma\mino_{\vu\in\Deltapl}
    [\partial_i F(\vu)]
\eeqn
$F_0=F(\vu^0)$, and
$\u^1_i=\delta_{ii^1}$ with $i^1=\argmax_i F_{iR}^{ub}$, and
$\u^2_i=\delta_{ii^2}$ with $i^2=\argmin_i F_{iR}^{lb}$, and $\leqsq$
defined in (\ref{eqleqsq}) means $\leq$ {\em and} $=\;+O(\sigma^2)$,
where $\sigma=1-\u_\p^0$.
\end{theorem}
For conservative estimates, the lower bound on $\low F$ and the upper
bound on $\up F$ are the interesting ones.

\paragraph{Proof} We start by giving an $O(\sigma^2)$ bound on
$\up F_R=\maxo_{\vu\in\Delta'}F_R(\vu)$. We first insert
(\ref{dtbnd}) with $\vt^*=\vt^0\equiv\v 0$ into
(\ref{Fexpand}) and treat $\v{\check\u}$ and $\vt$ as separate
variables:
\beqn
  F_R(\v{\check\u},\vt) \;=\; \sigma\sum_i
  [\partial_i F(\v{\check\u})]\cdot t_i
  \;\leqsq\; \maxo_{\v{\check\u}\in\Deltapl}
  \bigg\{\sigma\sum_i[\partial_i F(\v{\check\u})]\cdot t_i\bigg\}
  \;\leqsq\;
  \sum_i F_{iR}^{ub}\cdot t_i
\eeqn
\beq\label{defFiRub} %\label{F1gen}
  \qmbox{with}
  F_{iR}^{ub} \;:=\; \sigma\maxo_{\v{\check\u}\in\Deltapl}
    [\partial_i F(\v{\check\u})]
\eeq
The first inequality is obvious, the second follows from the
convexity of $\max$. From assumption (\ref{Lipschitz}) we get
$\partial_i F(\vu)-\partial_i F(\vu') = O(\sigma)$ for all
$\vu,\vu'\in\Deltapl$, since $\Deltapl$ has
dia\-meter $O(\sigma)$. Due to one additional $\sigma$ in
(\ref{defFiRub}) the expressions in (\ref{defFiRub}) change only by %{F1gen}
$O(\sigma^2)$ when introducing or dropping
$\maxo_{\v{\check\u}}$ anywhere. This shows that the
inequalities are tight within $O(\sigma^2)$ and justifies
$\leqsq$. We now upper bound $F_R(\vu)$:
\beq\label{FRup}
  \up F_R = \maxo_{\vu\in\Delta'} F_R(\vu)
  \leqsq \maxo_{\vt\in\Delta}\maxo_{\v{\check\u}\in\Deltapl}
  F_R(\v{\check\u},\vt)
  \leqsq
  \maxo_{\vt\in\Delta}\sum_i F_{iR}^{ub}\cdot t_i
  =
  \max_i F_{iR}^{ub} =: F_R^{ub}
\eeq
A linear function on $\Delta$ is maximized by setting the $t_i$
component with largest coefficient to 1. This shows the last
equality. The maximization over $\v{\check\u}$ in (\ref{defFiRub}) can
often be performed analytically, leaving an easy $O(d)$ time task
for maximizing over $i$.

We have derived an upper bound $F_R^{ub}$ on $\up F_R$. Let us
define the corner $t_i=\delta_{ii^1}$ of $\Delta$ with
$i^1:=\arg\max_i F_{iR}^{ub}$. Since $\up F_R\geq F_R(\vu)$
for all $\vu$, $F_R(\vu^1)$ in particular is a lower bound
on $\up F_R$. A similar line of reasoning as above shows that that
$F_R(\vu^1)=\up F_R+O(\sigma^2)$. Using $\up{F+const.}=\up
F+const.$ we get $O(\sigma^2)$ lower and upper bounds on $\up F$,
i.e.\ $F(\vu^1)\leqsq\up F\leqsq F_0+F_R^{ub}$. $\low F$ is bound
similarly with all max's replaced by min's and inequalities
reversed.
Together this proves the Theorem \ref{thARI}. \qed

%%%%%%%%%%%%%%%%%%%%%%%%%%%%%%%%%%%%%%%%%%%%%%%%%%%%%%%%%%%%%%%
\section{Error Propagation}\label{secEP}
%%%%%%%%%%%%%%%%%%%%%%%%%%%%%%%%%%%%%%%%%%%%%%%%%%%%%%%%%%%%%%%

%-------------------------------%
\subsection{Approximation of $\up{\low F}$ (special cases)}
%-------------------------------%
For the special case $F(\vu)=\sum_i f(\vu)$ we have $\partial_i
F(\vu)=f'(\u_i)$. For concave $f$ like in case of the entropy we
get particularly simple bounds
\beqn
  F_{iR}^{ub} = \sigma\maxo_{\vu\in\Deltapl} f'(\u_i)
%  = \sigma\max_{\u_i\in[\u_i^0,\u_i^0+\sigma]} f'(\u_i)
  = \sigma f'(\u_i^0),\quad\;\;
  F_R^{ub}=\sigma\max_i f'(\u_i^0)
  = \sigma f'(\textstyle{\min_i n_i\over n+s}),\quad\;\;
\eeqn
\beqn
  F_{iR}^{lb} = \sigma\mino_{\vu\in\Deltapl} f'(\u_i)
%  = \sigma\min_{\u_i[\u_i^0,\u_i^0+\sigma]} f'(\u_i)
  = \sigma f'(\u_i^0+\sigma),\quad
  F_R^{lb}=\sigma\min_i f'(\u_i^0+\sigma)
  = \sigma f'(\textstyle{\max_i n_i+s\over n+s}),
\eeqn
where we have used $\maxo_{\vu\in\Deltapl} f'(\u_i)
=\max_{\u_i\in[\u_i^0,\u_i^0+\sigma]} f'(\u_i)=f'(\u_i^0)$, and similarly
for $\min$. Analogous results hold for convex functions. In case
the maximum cannot be found exactly one is allowed to further
increase $\Deltapl$ as long as its diameter remains $O(\sigma)$.
Often an increase to
$\Box':=\{\vu:\u_i^0\leq\u_i\leq\u_i^0+\sigma\} \supset
\Deltapl \supset \Delta'$ makes the problem easy. Note that if we
were to perform these kind of crude enlargements on
$\maxo_{\vu}F(\vu)$ directly we would loose the bounds by
$O(\sigma)$.

\paragraph{Example (continued)}
$\sigma={1\over 10}$,
$h'({3\over 10})= {13051\over 2520}-\odt\pin^2$,
$h'({7\over 10})= {91717\over 8400}-{7\over 6}\pin^2$,
$H_0=H(\vu^0)=h({3\over 10})+h({6\over 10})$,
$H_R^{ub}={1\over 10}h'({3\over 10})$,
$H_R^{lb}={1\over 10}h'({7\over 10})$ $\Rightarrow$
$[H_0+H_R^{lb}\, ,\, H_0+H_R^{ub}]=[0.5564...,0.6404...]$, hence
$H_0+H_R^{ub}-\up H=0.0148=O({1\over 10^2})$,
$\low H-H_0-H_R^{lb}=0.0074...=O({1\over 10^2})$.

%-------------------------------%
\subsection{Error propagation}
%-------------------------------%
Assume we found bounds for estimators $G(\vu)$ and
$H(\vu)$ and we want now to bound the sum
$F(\vu):=G(\vu)+H(\vu)$. In the direct approach $\up
F\leq \up G+\up H$ we may lose $O(\sigma)$. A simple example is
$G(\vu)=\u_i$ and $H(\vu)=-\u_i$ for which
$F(\vu)=0$, hence $0=\up F\leq \up G+\up
H=\u_i^0+\sigma-\u_i^0=\sigma$, i.e.\ $\up F\not\leqsq \up G+\up
H$.
We can exploit the techniques of the previous section to obtain
$O(\sigma^2)$ approximations.
\beqn
  F_{iR}^{ub} \;=\;
  \sigma\maxo_{\vu\in\Deltapl}\partial_i F(\vu)
  \;\leqsq\;
  \sigma\maxo_{\vu\in\Deltapl}\partial_i G(\vu) +
  \sigma\maxo_{\vu\in\Deltapl}\partial_i H(\vu)
  \;=\; G_{iR}^{ub} + H_{iR}^{ub}
\eeqn
\begin{theorem}[Error propagation: Sum]\label{thEps}
Let $G(\vu)$ and $H(\vu)$ be Lipschitz differentiable and
$F(\vu)=\alpha G(\vu)+\beta H(\vu)$, $\alpha,\beta\geq 0$,
then $\up F\leqsq F_0+F_R^{ub}$ and $\low F\geqsq
F_0+F_R^{lb}$, where $F_0=\alpha G_0+\beta H_0$, and
$F_{iR}^{ub}\leqsq \alpha G_{iR}^{ub}+\beta H_{iR}^{ub}$, and
$F_{iR}^{lb}\geqsq \alpha G_{iR}^{lb}+\beta H_{iR}^{lb}$.
\end{theorem}
It is important to notice that $F_R^{ub}\not\leqsq
G_R^{ub}+H_R^{ub}$ (use previous example),
i.e.\ $\max_i[G_{iR}^{ub}+H_{iR}^{ub}]\not\leqsq\max_i
G_{iR}^{ub}+\max_i H_{iR}^{ub}$. $\max_i$ can not be pulled in and
it is important to propagate $F_{iR}^{ub}$, rather than
$F_R^{ub}$.

Every function $F$ with bounded curvature can be written as a sum
of a concave function $G$ and a convex function $H$. For convex and concave
functions, determining bounds is particularly easy, as we have seen.
Often $F$ decomposes naturally into convex and concave parts as is
the case for the mutual information, addressed later.
Bounds can also be derived for products.
\begin{theorem}[Error propagation: Product]\label{thEpp}
Let $G,H:\Deltapl\to[0,\infty)$ be non-nega\-tive
Lipschitz differentiable
functions (\ref{Lipschitz}) with non-negative
derivatives $\partial_i G,\partial_i H\geq 0$ $\forall i$ and
$F(\vu)=G(\vu)\cdot H(\vu)$, then $\up F\leqsq
F_0+F_R^{ub}$, where $F_0= G_0\cdot H_0$, and $F_{iR}^{ub}\leqsq
G_{iR}^{ub}(H_0+H_R^{ub})+ (G_0+G_R^{ub})H_{iR}^{ub}$, and
similarly for $\low F$.
\end{theorem}

\paragraph{Proof} We have
\beqn
  F_{iR}^{ub}
  \;=\; \sigma\max\partial_i F
  \;=\; \sigma\max\partial_i (G\!\cdot\!H)
  \;=\; \sigma\max[(\partial_i G)H+G(\partial_i H)]
  \;\leqsq\;
\eeqn
\beqn
         \sigma(\max\partial_i G)(\max H) +
         \sigma(\max G)(\max \partial_i H)
  \leqsq G_{iR}^{ub}(H_0\!+\!H_R^{ub}) + (G_0\!+\!G_R^{ub})H_{iR}^{ub}
\eeqn
where all functions depend on $\vu$ and all $\max$ are over
$\vu\in\Deltapl$. There is one subtlety in the last
inequality: $\max G\neq \up G\leqsq G_0+G_R^{ub}$.
The reason for the $\neq$ being that the maximization is taken
over $\Deltapl$, not over $\Delta'$ as in the definition
of $\up G$. The
correct line of reasoning is as follows:
\beqn
  \maxo_{\vu\in\Deltapl} G_R(\vu)
  \leqsq
  \maxo_{\vt\in\Deltal}\sum_i G_{iR}^{ub}\cdot t_i
  =
  \max\{0,\max_i G_{iR}^{ub}\} = G_R^{ub}
  \;\Rightarrow\;
  \max G\leqsq G_0+G_R^{ub}
\eeqn
The first inequality can be proven in the same way as
(\ref{FRup}). In the first equality we set the $t_i=1$ with
maximal $G_{iR}^{ub}$ {\em if} it is positive. If all
$G_{iR}^{ub}$ are negative we set $\vt\equiv\v 0$. We assumed
$G\geq 0$ and $\partial_i G\geq 0$, which implies $G_R\geq 0$. So,
since $G_R\geq 0$ anyway, this subtlety is ineffective. Similarly
for $\max H_R$.\qed

It is possible to remove the rather strong non-negativity
assumptions. Propagation of errors for other combinations like
ratios $F=G/H$ may also be obtained.

%%%%%%%%%%%%%%%%%%%%%%%%%%%%%%%%%%%%%%%%%%%%%%%%%%%%%%%%%%%%%%%
\section{Robust Intervals for Mutual Information}\label{secIEMI}
%%%%%%%%%%%%%%%%%%%%%%%%%%%%%%%%%%%%%%%%%%%%%%%%%%%%%%%%%%%%%%%

%-------------------------------%
\subsection{Mutual Information}
%-------------------------------%
We illustrate the application of the previous results on the
Mutual Information between two random variables
$\imath\in\{1,...,d_1\}$ and $\jmath\in \{1,...,d_2\}$. Consider
an i.i.d.\ random process with outcome
$(i,j)\in\{1,...,d_1\}\times\{1,...,d_2\}$ having joint
probability $\pi_{ij}$, where $\vpi\in\Deltapi :=\{\v x\in
I\!\!R^{d_1\times d_2}\,:\,x_{ij}\geq 0\,\forall ij,\; x_\pp=1\}$.
An important measure of the stochastic dependence of $\imath$ and
$\jmath$ is the mutual information
\beq\label{mi}
  {\cal I}(\vpi) = \sum_{i=1}^{d_1}\sum_{j=1}^{d_2}
  \pi_{ij}\log{\pi_{ij}\over\pi_{i\p}\pi_{\p j}} =
  \sum_{ij}\pi_{ij}\log\pi_{ij} -\!
  \sum_{i}\pi_{i\p}\log\pi_{i\p} -\!
  \sum_{j}\pi_{\p j}\log\pi_{\p j}
\eeq
\beqn
  \;=\; {\cal H}(\vpi_{\imath\p}) +
        {\cal H}(\vpi_{\p\jmath}) - {\cal
        H}(\vpi_{\imath\jmath})
\eeqn
$\pi_{i\p}=\sum_j\pi_{ij}$ and $\pi_{\p j}=\sum_i\pi_{ij}$ are row
and column marginal chances. Again, we assume a Dirichlet prior
over $\vpi_{\imath\jmath}$, which leads to a Dirichlet posterior
$p(\vpi_{\imath\jmath}|\v
n)\propto\prod_{ij}\pi_{ij}^{n_{ij}+st_{ij}-1}$ with $\vt\in
\Delta$. The expected value of $\pi_{ij}$ is
\beqn
  E_\vt[\pi_{ij}]={n_{ij}+st_{ij}\over n+s}=:\u_{ij}
\eeqn
The marginals $\vpi_{i\p}$ and $\vpi_{\p j}$ are also Dirichlet
with expectation $\u_{i\p}$ and $\u_{\p j}$.
The expected mutual information $I(\vu):=E_\vt[{\cal I}]$
can, hence, be expressed in terms of the expectations of
three entropies $H(\vu):=E_\vt[{\cal H}]$ (see (\ref{hex}))
\beqn
  I(\vu) = H(\vu_{\imath\p})+H(\vu_{\p\jmath})-H(\vu_{\imath\jmath})
 \;=\; H_{row}+H_{col}-H_{joint}
\eeqn
\beqn
 \;=\; \sum_i h(\u_{i\p})+\sum_j h(\u_{\p j})-\sum_{ij}h(\u_{ij})
\eeqn
where here and in the following we index quantities with $joint$,
$row$, and $col$ to denote to which distribution the quantity
refers.

%-------------------------------%
\subsection{Crude bounds for $I(\vu)$}
%-------------------------------%
Estimates for the robust IDM interval $[\mino_{\vt\in\Delta}E_\vt[{\cal I}] \,,\, \maxo_{\vt\in\Delta }E_\vt[{\cal I}]]$ can
be obtained by [minimizing,maximizing] $I(\vu)$. A crude upper
bound can be obtained as
\beqn
  \up I \;:=\; \maxo_{\vt\in\Delta} I(\vu) \;=\;
  \max[H_{row}+H_{col}-H_{joint}] \;\leq\;
\eeqn
\beqn
  \max H_{row} + \max H_{col} - \min H_{joint} \;=\;
  \up H_{row} + \up H_{col} -  \low H_{joint},
\eeqn
where exact solutions to $\up H_{row}$, $\low H_{row}$ and
$\low H_{joint}$ are available from Section \ref{secEEI}.
Similarly
$
  \low I \geq \low H_{row} + \low H_{col} -  \up H_{joint}
$. The problem with these bounds is that, although good in some
cases, they can become arbitrarily crude.
The following $O(\sigma^2)$ bound can be derived by exploiting the
error sum propagation Theorem \ref{thEps}.

\begin{theorem}[Bound on lower and upper Mutual Information]\label{thMIbnd}
The following bounds on the expected mutual information
$I(\vu)=E_\vt[{\cal I}]$ are valid:
\bqan
 & & \nq I(\vu^1) \leqsq \up I\leqsq I_0+I_R^{ub} \qmbox{and}
 I_0+I_R^{lb} \leqsq \low I \leqsq  I(\vu^2),
 \qmbox{where} \\
 & & \nq I_0=I(\vu^0)=H_{0row}+H_{0col}-H_{0joint} =
 h(\u_{i\p}^0)+h(\u_{\p j}^0)-h(\u_{ij}^0), \\
 & & \nq I_{ijR}^{ub}\leqsq
 H_{iRrow}^{ub}+H_{jRcol}^{ub}-H_{ijRjoint}^{lb} =
 h'(\u_{i\p}^0)+h'(\u_{\p j}^0)-h'(\u_{ij}^0\!+\!\sigma), \\
 & & \nq I_{ijR}^{lb}\geqsq
 H_{iRrow}^{lb}+H_{jRcol}^{lb}-H_{ijRjoint}^{ub} =
 h'(\u_{i\p}^0\!+\!\sigma)+h'(\u_{\p j}^0\!+\!\sigma)-h'(\u_{ij}^0),
\eqan
with $h$ defined in (\ref{hex}), and $t^0_{ij}=0$, and
$t^1_{ij}=\delta_{(ij)(ij)^1}$ with $(ij)^1=\argmax_{ij}
I_{ijR}^{ub}$, and $t^2_{ij}=\delta_{(ij)(ij)^2}$ with
$(ij)^2=\argmin_{ij} I_{ijR}^{lb}$.
\end{theorem}

%%%%%%%%%%%%%%%%%%%%%%%%%%%%%%%%%%%%%%%%%%%%%%%%%%%%%%%%%%%%%%%
\section{IDM for Product Spaces}\label{secPS}
%%%%%%%%%%%%%%%%%%%%%%%%%%%%%%%%%%%%%%%%%%%%%%%%%%%%%%%%%%%%%%%

%-------------------------------%
%\subsection{IDM for product spaces}
%-------------------------------%
Product spaces $\Omega=\Omega_1\times...\times\Omega_m$ with
$\Omega_k=\{1,...d_k\}$ occur frequently in practical problems,
e.g.\ in the mutual information ($m=2$), in robust trees ($m=3$),
or in Bayesian nets in general ($m$ large). Without loss of
generality we only discuss the $m=2$ case in the following.
Ignoring the underlying structure in $\Omega$, a Dirichlet prior
in case of unknown chances $\pi_{\imath\jmath}$ and an IDM as used
in Section \ref{secIEMI} with
\beq\label{IDMfull}
\vt\in\Delta :=\{\vt\in I\!\!R^{d_1\times
d_2}\equiv I\!\!R^{d_1}\otimes I\!\!R^{d_2}\,:\,t_{ij}\geq
0\,\forall ij,\; t_\pp=1\}
\eeq
seems natural. On the other hand, if we take into account the
structure of $\Omega$ and go back to the original motivation of
IDM this choice is far less obvious. Recall that one of the major
motivations of IDM was its reparametrization invariance in the
sense that inferences are not affected when grouping or splitting
events in $\Omega$. For unstructured spaces like $\Omega_k$ this
is a reasonable principle. For illustration, let us consider
objects of various {\em shape} and {\em color}, i.e.\
$\Omega=\Omega_1\times\Omega_2$, $\Omega_1=\{ball, pen, die,
...\}$, $\Omega_2=\{yellow, red, green, ...\}$ in generalization
to Walleys bag of marbles example. Assume we want to detect a
potential dependency between {\em shape} and {\em color} by means
of their mutual information $I$. If we have no prior idea on the
possible kind of colors, a model which is independent of the
choice of $\Omega_2$ is welcome. Grouping red and green, for
instance, corresponds to $(x_{i1}$, $x_{i2}$, $x_{i3}$, $x_{i4},
...)$ $\leadsto$ $(x_{i1}$, $x_{i2}+x_{i3}$, $x_{i4}, ...)$ {\em
for all shapes $i$}, where $\v x\in\{\v n,\vpi,\vt,\vu\}$.
Similarly for the different shapes, for instance we could group
all round or all angular objects. The ``smallest IDM'' which
respects this invariance is the one which considers all
\beq\label{IDMprod}
  \vt\in\DeltaOX:=\Delta_{d_1}\otimes\Delta_{d_2}
  \;\propersubset\;\Delta.
\eeq
The tensor or outer product $\otimes$ is defined as $(\v
v\otimes\v w)_{ij}:=v_iw_j$ and $V\otimes W:=\{\v v\otimes\v w:\v
v\in V,\, \v w\in W\}$. It is a bilinear (not linear!) mapping.
This ``small tensor'' IDM is invariant under arbitrary grouping of
columns and rows of the chance matrix $(\vpi_{ij})_{1\leq i\leq
d_1,1\leq j\leq d_2}$. In contrast to the larger $\Delta$ IDM
model it is not invariant under arbitrary grouping of matrix
cells, but there is anyway little motivation for the necessity of
such a general invariance. General non-column/row cross groupings
would destroy the product structure of $\Omega$ and with that the
mere concepts of shape and color, and their correlation. For $m>2$
as in Bayes-nets cross groupings look even less natural. Whether
the $\DeltaOX$ or the larger simplex $\Delta$ is the more
appropriate IDM model depends on whether one regards the structure
$\Omega_1\times\Omega_2$ of $\Omega$ as a natural prior knowledge
or as an arbitrary a posteriori choice. The smaller IDM has the
potential advantage of leading to more precise
predictions (smaller robust sets).

Let us consider an estimator $F:\Delta\to\Set R$ and its
restriction $F_\ots:\DeltaOX\to\Set R$.
Robust intervals $[\low F,\up F]$ for $\Delta$ are generally wider
than robust intervals $[\low{F}_\ots,\up F_\ots]$ for $\DeltaOX$.
Fortunately not much. Although $\DeltaOX$ is a {\em
lower-dimensional} subspace of $\Delta$, it contains all vertices
of $\Delta$. This is possible since $\DeltaOX$ is a {\em
nonlinear} subspace. The set of ``vertices'' in both cases is
$\{\vt\,:\,t_{ij}=\delta_{ii_0}\delta_{jj_0},\;i_0\in\Omega_1,\;j_0\in\Omega_2\}$.
Hence, {\em if} the robust interval boundaries $\up{\low F}$ are
assumed in the vertices of $\Delta$ {\em then} the
interval for the $\DeltaOX$ IDM model is the
same ($\up{\low F}=\up{\low F}_\ots$). Since the condition is
``approximately'' true, the conclusion is
``approximately'' true. More precisely:

\begin{theorem}[IDM bounds for product spaces]
The $O(\sigma^2)$ bounds of Theorem \ref{thARI} on the robust
interval $\up{\low F}$ in the full IDM model $\Delta$
(\ref{IDMfull}), remain valid for $\up{\low
F}_\ots$ in the product IDM model
$\DeltaOX$
(\ref{IDMprod}).
\end{theorem}

\paragraph{Proof}\vspace{-3ex}
\beqn
  F(\vu^1) \leq \up F_\ots \leq \up F \leq
  F_0+F_R^{ub} = F(\vu^1)+O(\sigma^2),
\eeqn
where $\up F_\ots:=\maxo_{\vt\in\deltaOX}F(\vu)$ and $\vu^1$
was the ``$F_R$ maximizing'' vertex as defined in Theorem
\ref{thMIbnd} ($F(\vu^1)\leqsq\up F$). The first inequality
follows from the fact that all $\Delta$ vertices also belong to
$\DeltaOX$, i.e.\ $\vt^1\in\DeltaOX$. The second inequality
follows from $\DeltaOX\subset\Delta$. The remaining (in)equalities
follow from Theorem \ref{thARI}. This shows that $|\up F_\ots-\up
F|=O(\sigma^2)$, hence $F_0+F_R^{ub}$ is also an $O(\sigma^2)$
upper bound to $\up F_\ots$. This implies that to the
approximation accuracy we can achieve, the choice between $\Delta$
and $\DeltaOX$ is irrelevant. \qed

%%%%%%%%%%%%%%%%%%%%%%%%%%%%%%%%%%%%%%%%%%%%%%%%%%%%%%%%%%%%%%%
\section{Robust Credible Intervals}\label{secCI}
%%%%%%%%%%%%%%%%%%%%%%%%%%%%%%%%%%%%%%%%%%%%%%%%%%%%%%%%%%%%%%%

%-------------------------------%
\subsection{Bayesian credible sets/intervals}
%-------------------------------%
For a probability distribution $p:I\!\!R^d\to[0,1]$, an
$\alpha$-credible region is a measurable set $A$ for which
$p(A):=\int p(x)\chi_A(x) d^dx\geq\alpha$, where $\chi_A(x)=1$ if
$x\in A$ and $0$ otherwise, i.e.\ $x\in A$ with probability at
least $\alpha$. For given $\alpha$, there are many choices for
$A$. Often one is interested in ``small'' sets, where the size of
$A$ may be measured by its volume $\Vol(A):=\int\chi_A(x)d^dx$.
Let us define a/the smallest $\alpha$-credible set
\beqn
  A^{min}:=\argmin_{A:p(A)\geq\alpha}\Vol(A)
\eeqn
with ties broken arbitrarily. For unimodal $p$, $A^{min}$ can be
chosen as a connected set. For $d=1$ this means that
$A^{min}=[a,b]$ with $\int_a^b p(x)dx=\alpha$ is a minimal length
$\alpha$-credible interval. If, additionally $p$ is symmetric
around $E[x]$, then $A^{min}=[E[x]-a,E[x]+a]$ is also symmetric
around $E[x]$.

%-------------------------------%
\subsection{Robust credible sets}
%-------------------------------%
If we have a set of probability distributions $\{p_t(x)$, $t\in
T\}$, we can choose for each $t$ an $\alpha$-credible set $A_t$
with $p_t(A_t)\geq\alpha$, a minimal one being
$A_t^{min}:=\argmin_{A:p_t(A)\geq\alpha}\Vol(A)$. A robust
$\alpha$-credible set is a set $A$ which contains $x$ with
$p_t$-probability at least $\alpha$ for {\em all} $t$. A minimal
size robust $\alpha$-credible set is
\beq\label{eqRCI}
  A^{min}:=\argmin_{A={\textstyle\cup}_t A_t:p_t(A_t)\geq\alpha\forall t\in T}
  \Vol(A)
\eeq
It is not easy to deal with this expression, since $A^{min}$ is
{\em not} a function of $\{A_t^{min}:t\in T\}$, and especially does
not coincide with $\bigcup_t A_t^{min}$ as one might expect.

%-------------------------------%
\subsection{Robust credible intervals}
%-------------------------------%
This can most easily be seen for univariate symmetric
unimodal distributions, where $t$ is a translation, e.g.\
$p_t(x)=\mbox{Normal}(E_t[x]=t,\sigma=1)$ with 95\% credible
intervals $A_t^{min}=[t-2,t+2]$. For, e.g.\ $T=[-1,1]$ we get
$\bigcup_t A_t^{min}=[-3,3]$. The credible intervals {\em move}
with $t$. One can get a smaller union if we take the intervals
$A_t^{sym}=[-a_t,a_t]$ symmetric around 0. Since $A_t^{sym}$ is a
non-central interval w.r.t.\ $p_t$ for $t\neq 0$, we have $a_t>2$, i.e.\
$A_t^{sym}$ is larger than $A_t^{min}$, but one can show that the
increase of $a_t$ is smaller than the shift of $A_t^{min}$ by $t$,
hence we save something in the union. The optimal choice is
neither $A_t^{sym}$ nor $A_t^{min}$, but something in-between.
In the extended version \cite{Hutter:03idmx} this
is illustrated for the triangular distribution
$p_t(x)=\max\{0\,,\,1\!-\!|x\!-\!t|\}$ with $t\in T:=[-\gamma,\gamma]$,
where closed form solutions can be given.

An interesting open question is under which general conditions we
can expect $A^{min}\subseteq\bigcup_t A_t^{min}$. In any case,
$\bigcup_t A_t$ can be used as a conservative estimate for a
robust credible set, since $p_t(\bigcup_{t'} A_{t'})\geq
p_t(A_t) \geq \alpha$ for all $t$.

A special (but important) case which falls outside the above
framework are one-sided credible intervals, where only $A_t$ of
the form $[a,\infty)$ are considered. In this case
$A^{min}=\bigcup_t A_t^{min}$, i.e.\ $A^{min}=[a_{min},\infty)$
with $a_{min}=\max\{a:p_t([a,\infty])\geq\alpha\forall t\}$.

%-------------------------------%
\subsection{Approximations}
%-------------------------------%
For complex distributions like for the mutual information we have
to approximate (\ref{eqRCI}) somehow. We use the following notation
for shortest $\alpha$-credible {\em intervals} w.r.t.\
a univariate distribution $p_t(x)$:
\beqn
  {\mathop{\widetilde{x}}\limits_\sim}\!\,_t \;\equiv\;
  [ {\mathop{x}\limits_\sim}\!\,_t , \widetilde{x}_t] \;\equiv\;
  [ E_t[x]-\Delta{\mathop{x}\limits_\sim}\!\,_t \,,\,
  E_t[x]+\Delta\widetilde{x}_t ] \;:=\;
  \argmin_{[a,b]:p_t([a,b])\geq\alpha}(b-a),
\eeqn
where $\Delta\widetilde{x}_t:=\widetilde{x}_t-E_t[x]$
($\Delta{\mathop{x}\limits_\sim}\!\,_t :=
E_t[x]-{\mathop{x}\limits_\sim}\!\,_t$) is the distance from the
right boundary $\widetilde{x}_t$ (left boundary
${\mathop{x}\limits_\sim}\!\,_t$) of the shortest
$\alpha$-credible interval
${\mathop{\widetilde{x}}\limits_\sim}\!\,_t$ to the mean $E_t[x]$
of distribution $p_t$. We can use
$\mathop{\up{\widetilde{x}}}\limits_\simeq \equiv
[\mathop{x}\limits_\simeq , \up{\widetilde{x}}] :=\bigcup_t
{\mathop{\widetilde{x}}\limits_\sim}\!\,_t$ as a (conservative,
but not shortest) robust credible interval, since
$p_t(\mathop{\up{\widetilde{x}}}\limits_\simeq)\geq p_t
({\mathop{\widetilde{x}}\limits_\sim}\!\,_t)\geq\alpha$ for all
$t$. We can upper bound $\up{\widetilde{x}}$ (and similarly lower
bound $\mathop{x}\limits_\simeq$) by
\beq\label{eqCRCI}
  \up{\widetilde{x}} \;=\;
  \max_t(E_t[x]+\Delta\widetilde{x}_t) \;\leq\;
  \max_t E_t[x]+\max_t \Delta\widetilde{x}_t \;=\;
  \up{E[x]}+ \up{\Delta\widetilde{x}}.
\eeq
We have already intensively discussed how to compute upper and
lower quantities, particularly for the upper mean $\up{E[x]}$ for
$x\in\{{\cal F}, {\cal H}, {\cal I}, ... \}$, but the linearization
technique introduced in Section \ref{secAEI} is general enough to
deal with all in $t$ differentiable quantities, including
$\Delta\widetilde{x}_t$. For example for Gaussian $p_t$ with
variances $\sigma_t$ we have
$\Delta\widetilde{x}_t=\kappa\sigma_t$ with $\kappa$ given by
$\alpha=\mbox{erf}(\kappa/\sqrt{2})$, where erf is the error
function (e.g.\ $\kappa=2$ for $\alpha\approx 95\%)$. We only need to
estimate $\max_t\sigma_t$.

For non-Gaussian distributions, exact expression for
$\Delta\widetilde{x}_t$ are often hard or impossible to obtain and
to deal with. Non-Gaussian distributions depending on some sample
size $n$ are usually close to Gaussian for large $n$ due to the
central limit theorem.
One may simply use $\kappa\sigma_t$ in place of
$\Delta\widetilde{x}_t$ also in this case, keeping in mind that
this could be a non-conservative approximation. More
systematically, simple (and for large $n$ good) upper bounds on
$\Delta\widetilde{x}_t$ can often be obtained and should preferably be used.

Further, we have seen that the variation of sample depending
differentiable functions (like $E_t[x]=E_t[x|\v n]$) w.r.t.\ $t\in\Delta$
are of order ${s\over n+s}$. Since in such cases the standard
deviation $\sigma_t\sim n^{-1/2} \sim \Delta\widetilde{x}_t$ is
itself suppressed, the variation of $\Delta\widetilde{x}_t$ with
$t$ is of order $n^{-3/2}$. If we regard this as negligibly
small, we may simply fix some $t^*\in\Delta$:
\beqn
  \max_t\Delta\widetilde{x}_t = \kappa\sigma_{t^*}+O(n^{-3/2})
\eeqn
Since $\Delta\widetilde{x}_t$ is ``nearly'' constant, this also
shows that we lose at most $O(n^{-3/2})$ precision in the bound
(\ref{eqCRCI}) (equality holds for $\Delta\widetilde{x}_t$
independent of $t$).
Expressions for the variance of $\cal I$, for instance, have been
derived in \cite{Wolpert:95,Hutter:01xentropy}:

%%%%%%%%%%%%%%%%%%%%%%%%%%%%%%%%%%%%%%%%%%%%%%%%%%%%%%%%%%%%%%%
\section{Conclusions}\label{secConc}
%%%%%%%%%%%%%%%%%%%%%%%%%%%%%%%%%%%%%%%%%%%%%%%%%%%%%%%%%%%%%%%

% Section 1 & 2
This is the first work, providing a systematic approach for
deriving closed form expressions for interval estimates in the
Imprecise Dirichlet Model (IDM).
% Section 3 & 4
We concentrated on exact and conservative {\em robust} interval
([lower,upper]) estimates for concave functions
$F=\sum_i f_i$ on simplices, like the entropy. The
conservative estimates widened the intervals by $O(n^{-2})$, where
$n$ is the sample size.
Here is a dilemma, of course: For large $n$ the approximations are
good, whereas for small $n$ the bounds are more interesting, so
the approximations will be most useful for intermediate $n$. More
precise expressions for small $n$ would be highly interesting.
% Section 5 and 6
We have also indicated how to
propagate robust estimates from simple functions to composite
functions, like the mutual information.
% Section 7
We argued that a reduced IDM on product spaces, like Bayesian
nets, is more natural and should be preferred in order to improve
predictions. Although improvement is formally only $O(n^{-2})$,
the difference may be significant in Bayes nets or for very small $n$.
Finally, the basics of how to combine robust with credible
intervals have been laid out. Under certain conditions
$O(n^{-3/2})$ approximations can be derived, but the presented
approximations are not conservative.
All in all this work has shown that IDM has not only interesting
theoretical properties, but that explicit
(exact/conservative/approximate) expressions for robust
(credible) intervals for various quantities can be derived. The
computational complexity of the derived bounds on $F=\sum_i f_i$
is very small, typically one or two evaluations of $F$ or related
functions, like its derivative.
First applications of these (or more precisely, very similar)
results, especially the mutual information, to robust
inference of trees look promising \cite{Hutter:03tree}.

%------------------------------%
\subsection{Acknowledgements}
%------------------------------%
I want to thank Peter Walley for introducing IDM to me and
Marco Zaffalon for encouraging me to investigate this topic.
This work was supported by SNF grant 2000-61847.00 to
J{\"u}rgen Schmidhuber.

%%%%%%%%%%%%%%%%%%%%%%%%%%%%%%%%%%%%%%%%%%%%%%%%%%%%%%%%%%%%%%%
%         Bibliography        %
%%%%%%%%%%%%%%%%%%%%%%%%%%%%%%%%%%%%%%%%%%%%%%%%%%%%%%%%%%%%%%%

\end{document}